\documentclass{article}
\usepackage{arxiv}
\usepackage[utf8]{inputenc}
\usepackage[T1]{fontenc} 
\usepackage{hyperref}   
\usepackage{url}   
\usepackage{booktabs}
\usepackage{amsfonts}
\usepackage{nicefrac}
\usepackage{microtype}
\usepackage{lipsum}
\usepackage{graphicx}
\graphicspath{ {./images/} }

\title{
Synchronous Heterogeneous Exclusion Processes on Open Lattice
}

\author{
Marina V. Yashina  \\
  Department of Higher Mathematics\\
  Moscow Automobile and Road Construction\\
  State Technical University (MADI) \\
  Moscow, Leningradsky avenue, 64, Russia  \\
  \texttt{mv.yashina@madi.ru} \\ 
\And
 Alexander G. Tatashev  \\
  Department of Higher Mathematics\\
  Moscow Automobile and Road Construction\\
  State Technical University (MADI) \\
  Moscow, Leningradsky avenue, 64, Russia  \\
  \texttt{a-tatashev@yandex.ru} \\
}

\begin{document}
\maketitle
\begin{abstract}
 
A traffic model on an open one-dimensional lattice is considered. At any discrete time moment, with prescribed probability, a particle arrives to the leftmost cell of the lattice, and, with  prescribed probability, the arriving particle belongs to one of the types characterized by the probabilities of particle attempts to move at the present time and the probabilities to leave the system. An approximate approach to compute the particle flow rate and density in cells is proposed. It is proven that, for a particular case of the system, the approach gives exact results.

\end{abstract}

\keywords{ Traffic flow models \and Dynamical systems \and Synchronous exclusion processes \and  Spectrum of Buslaev net  \and  Algebraic structures
}

\section*{1.Introduction}
In the well-known Nagel-Schreckenberg traffic model [1], particles corresponding to vehicles move on a lattice according to prescribed rules. In any cell of the lattice, no more than one particle can be located simultaneously. The model can be interpreted as a cellular automaton [2] or as a stochastic exclusion process [3]. In the version of the random process introduced in [3], the time scale is continuous, and, in a time interval of duration $\Delta t,$ a particle tries to move with probability equal to $\lambda \Delta t+o(\Delta t),$ $\Delta t\to 0.$ This is equivalent to that the duration of the time interval between particle attempts to move is distributed exponentially with the average value $1/\lambda.$ This exclusion processes with continuous time is called asynchronous [4]. For an asynchronous exclusion  process, the probability of simultaneous movement of more than one particle is zero. Another version of exclusion process is called a synchronous exclusion process [4]. For this version, the time scale is discrete and it is possible that more than one particle move simultaneously. 

In [5], a simple version of the Nagel--Schreckenberg model is introduced. In the model, $M$ particles are located on a one-dimensional closed lattice containing $N$ cells. In this 
synchronous version, at any discrete time moment, a particle moves one cell forward with a prescribed probability $p$ under the assumption that the cell ahead is vacant.  The study of analogous asynchronous model is trivial, since, by considering the corresponding Markov chain, it is easy to verify that the stationary probabilities for all particle configurations are the same. As it is noted in [5], a similar property of equiprobability of particle configurations is possessed by a discrete-time model in that, at any discrete moment, only one particle moves with prescribed probability to a vacant cell, and any particle is chosen equiprobably (another version of an asynchronous exclusion process).  The average number of a particle transitions per unit time, which for the synchronous model, is equal to the stationary probability of a particle movement at a discrete moment is called the average velocity of particles. In [5], for the synchronous model, a formula for the average velocity was obtained, based on intuitive reasoning. As it is proved in [6], this formula gives the value of the average velocity on the related infinite lattice. In [6], for the synchronous model introduced in [5], an algorithm was found that allows to compute the value of the average particle velocity for prescribed $N$ and $M.$ In [7], an explicit formula for the value of the average velocity was obtained. The model introduced in [5] is stochastically equivalent to the zero-range process studied in [8]. In [9], [10], the behavior of the model introduced in [5] was considered in the limiting case $p=1.$ This deterministic version of the model corresponds to an elementary cellular automaton (ECA) rule~184 in the S.~Wolfram classification (ECA~184). Some versions and generalizations of the model introduced in [5] were studied in [4]. In [11], a generalization of the synchronous model introduced in [5] was considered for the case in that there are several types of particles characterized by different probabilities of moving to a vacant cell, as well as the corresponding model with continuous time (asynchronous model) in that the types of particles differ in the rate of particle attempts to move. In [11], formulas for the average particle velocity are obtained, and known results for queueing networks with discrete [12] and continuous [13], [14] time were used to derive these formulas. In [15], a generalization was obtained for the asynchronous model considered in [11]. In the model studied in [15], the particles can move in both directions.

In [16], a model was considered in that partiles move on two closed sequences of particles (rings) according to the rule of ECA~184. The rings can differ in length (number of cells). There is a common cell of the rings~--- a node. Passing through the node, a particle can move from one ring to another.

In [17], a model on one-dimensional open lattice is considered. The lattice contains a prescribed number of cells. The time is discrete. At each step, a particle enters the leftmost cell with a given probability under the assumption that this cell is vacant. At any step, each particle moves to the cell  adjecent to the right under the assumption  that this cell is vacant. At the current discrete moment, the particle leaves the system from the rightmost cell with prescribed probability. In [18], an asynchronous version of the system considered in [17] was studied. In [19], a two-lane version of the model [18] was considered.

In this paper, we consider a generalization of the model considered in [17] to the case of different types of particles (heterogeneous particle flow). 
An approach has been developed that allows to find the stationary probability of occupancy (density) for each prescribed cell and the particle flow rate.

In Section 2, an information is provided information on a synchronous model on an open lattice with one type of particle.

In Section 3, the considered synchronous model on an open lattice with several types of particles is described. 

In Section 4, the considered system is presented as a Markov chain and the ergodicity of the system is proved.

In Section 5, a system with one type of particles is introduced. This system corresponds to the original system with different types of particles. The introduced system is used as an auxiliary system in the study of the original system.

In Section 6, for a special case of the considered system, theorems
on the characteristics of the system are proved.

In Section 7, an approximate approach is proposed for calculating the particle density depending on the location on the lattice and the particle flow rate in the model. The approximate values obtained with aid of the proposed  approach are compared with the exact values obtained from  stationary state probabilities equations for the Markov chain corresponding to the considered models. For the special case considered in Section~6, the values of the system characteristics obtained with aid of the proposed approach are equal to the exact values.

\section*{2.
Information on the synchronous exclusion process on an open lattice
}
\label{section:BuslN}

\hskip 16pt The following system was studied in [17]. There is an open lattice containing $N$ cells with indices $1,\dots,N.$ Time is discrete. If at time $t=0,1,2,\dots,$ the cell 1 is vacant, then,  at time $t+1,$ with probability $\alpha,$ there will be a particle arriving from outside to this cell.  If, at time $t,$ a particle is in the cell $i,$ and the cell $i+1$ is vacant, then, at time $t+1,$ with probability~$p,$ the particle will be in cell~$i+1,$ $i=1,\dots,N-1.$ If, at time~$t,$ the particle is in cell $N,$ then, with probability~$\beta,$ the particle leaves the system and, at time~$t+1,$ the cell $N$ will be vacant.

In [17], an algorithm was developed that allows to compute the stationary probabilities of the 
system states, the stationary probabilities of occupancy (density) of cells depending on their indices, and the rate of the particle flow.

\section*{3. Description of synchronous exclusion processes with different types of particles}

\hskip 16pt Let us describe a generalization of the system described in 
Section~2. In the system, described in this section, there are several types of 
particles.  

The lattice of the system contains $N$ cells. The indices of the cells are $1,\dots,N.$ If, at time $t=0,1,2,\dots,$ the cell~1 is vacant, then, with probability $\alpha,$ at time $t+1,$ there will be a 
particle arriving from outside. There are $K$ types of particles. The probability that 
the arriving particle belongs to the type $k$ is $a_k,$ $k=1,\dots,K,$
$a_1+\dots+a_k=1.$ If, at time $t,$ the particle of the type $k$ is in the cell $i,$ and the cell $i+1$ is vacant, then, with probability~$p_k,$ at time $t+1,$ the particle will be
in the cell~$i+1,$ $k=1,\dots,K,$ $i=1,\dots,N-1.$ If, at time~$t,$ a particle of the type~$k$ is in the cell~$N,$ then,  with the probability~$\beta_k,$ the particle leaves the system, and, at time~$t+1,$ the cell $N$ is vacant, $k=1,\dots,K.$

\section*{4. Stationary state probabilities and ergodicity of system}

\hskip 16pt Let $x=(x_1,\dots,x_N)$ be a state of the system such that, if the cell~$i$ is vacant, 
then $x_i=0,$ and, if there is a particle of the type $k$ in the cell $i,$ then $x_i=k,$
$i=1,\dots,N,$ $k=1,\dots,K$, Figure 1. Let $x(t)$ be the state of the system at time $t\ge 0.$ The stochastic process $x(t)$ is a Markov chain with discrete time.

\begin{figure}[ht!]
\centerline{\includegraphics[width=300pt]{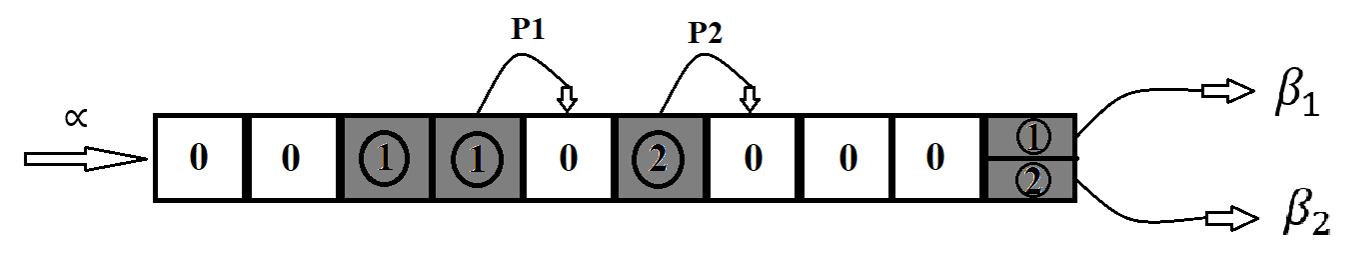}}
\caption{State (0,0,1,1,0,2,0,0,0,1|2), N=10, K=2}
\end{figure}

\vskip 3pt
{\bf Theorem 1.} {\it The system under consideration is ergodic, i.~e., for any state $(x_1,\dots,x_N),$ there exists the limit of the probability that the system is in this state as $t\to \infty$ (the stationary probability of the state), and this limit does not depend on the initial state.
\vskip 3pt
Proof.} The number of states of this chain is finite (this number is equal to $M=(k+1)^N),$ and the states form a unique non-periodic class of communicating states [20]. Indeed, from any state, the system results in the state $(0,0,\dots,0)$ after time interval with finite expectation, and, from the state $(0,0,\dots,0),$ the system results in any other state after time interval with finite expectation. This class of states is non-periodic because, at each step, the system can remain in the same state with positive probability. According to the ergodic theorem [20], for a finite Markov chain with a unique non-periodic class of communicating states, the following holds. For any state $(x_1,\dots,x_N),$ there exists a limit of the probability that, at time~$t,$ the system is in this state  as $t\to \infty,$ and this limit does not depend on the initial state. The theorem has been proved.
\vskip 3pt
Let the states of the chain be arbitrarily numbered. Let $P_l$ be the stationary probability of 
the state $l,$ $l=1,2,\dots,M.$ Let $p_{lr}$ be the probability of transition from the state $l$ to the state $r$ in one step, $l,r=1,2,\dots,M.$

According to the theory of Markov chains, stationary probabilities satisfy the system of linear equations
$$P_l\sum\limits_{r\ne l}p_{lr}=\sum\limits_{r\ne l}P_r p_{rl},\ l=1,\dots,M,$$
$$\sum\limits_{l=1}^MP_l=1.$$

Let the considered model be called the system $S.$

\section*{5. Description of auxiliary system}

\hskip 16pt We introduce the system $S^*$ that differs from the system $S$ in the 
following. In the system~$S^*,$ there is a unique type of particles, and the probability of a particle attempt to move, which we denote by~$p^*,$ is
computed according to the formula
$$p^*=\frac{1}{\sum\limits_{k=1}^K\frac{a_k}{p_k}},\eqno(1)$$
i.~e. $p^*$ is the harmonic average of the values $p_1,\dots,p_K.$ 

Similarly, suppose
$$\beta^*=\frac{1}{\sum\limits_{k=1}^K\frac{a_k}{\beta_k}}.\eqno(2)$$

The Markov chain corresponding to the system $S^*$ contains $M^*=2^N$ states.

Let $\eta(k)$ be a function defined on the set of numbers $0,1,\dots,K,$
$$\eta(k)=\left\{
\begin{array}{l}
0,\ k=0,\\
1,\ k>0.\\
\end{array}
\right.
$$

Suppose the state $(x_1^*,\dots,x_N^*)$ of the system $S^*$ is the state such that $x_i^*=0$ if the cell $i$ is vacant, and $x_i^*=1$ if the cell $i$ is occupied, $i=1,\dots,N.$

Suppose the state $(x_1^*,\dots,x_N^*)$ of the system $S^*$ corresponds to the set of the system $S$ states $(x_1,\dots,x_N)$ such that $(\eta(x_1),\dots,\eta(x_N))=(x_1^*,\dots,x_N^*).$ This set is called the set $G(x_1^*,\dots,x_N^*).$

\section*{6. Exact solution for a special case}

\hskip 16pt We will prove a theorem that allows us to reduce, in a special case, the computation of the stationary probabilities of cell occupancy and particle flow rate for the considered system $S$ with several types of particles to the computation of the corresponding characteristics of the system $S^*$ defined in Section~5.

Suppose $N=2$ and $\beta_1=\beta_2=\beta.$  Let $P(x_1,x_2)$ be the stationary probability that the system $S$ is in the state $(x_1,x_2).$ Let $P^*(x_1^*,x_2^*)$ be the stationary probability that the system $S$ is in the state $(x_1^*,x_2^*).$
\vskip 3pt
{\bf Theorem 2.} {\it The stationary probabilities  of the system $S$ states are expressed through the stationary probabilities of the the system $S^*$ states by equations
$$P(0,0)=P^*(0,0),\eqno(3)$$
$$P(0,1)=a_1P^*(0,1),\eqno(4)$$
$$P(0,2)=a_2P^*(0,1),\eqno(5)$$
$$P(1,0)=\frac{a_1p^*}{p_1}P^*(1,0),\eqno(6)$$
$$P(1,1)=a_1^2P^*(1,1),\eqno(7)$$
$$P(1,2)=a_1a_2P^*(1,1),\eqno(8)$$
$$P(2,0)=\frac{a_2p^*}{p_2}P^*(1,0),\eqno(9)$$
$$P(2,1)=a_1a_2P^*(1,1),\eqno(10)$$
$$P(2,2)=a_2^2P^*(1,1).\eqno(11)$$
\vskip 3pt
Proof.} Stationary states probabilities of the system $S^*$ satisfy the system of equations 
$$\alpha P^*(0,0)=(1-\alpha)\beta P^*(0,1),\eqno(12)$$ 
$$(1-(1-\alpha)(1-\beta))P^*(0,1)=p^* P^*(1,0),\eqno(13)$$
$$p^*P^*(1,0)=\alpha P^*(0,0)+\alpha \beta P^*(0,1)+\beta P^*(1,1),\eqno(14)$$
$$\beta P^*(1,1)=\alpha(1-\beta)P^*(0,1),\eqno(15)$$
$$P^*(0,0)+P^*(0,1)+P^*(1,0)+P^*(1,1)=1.\eqno(16)$$

Stationary states of the system $S$ satisfy the system of equations 
$$\alpha P(0,0)=(1-\alpha)\beta P(0,1)+(1-\alpha)\beta P(0,2),\eqno(17)$$ 
$$(1-(1-\alpha)(1-\beta))P(0,1)=p_1P(1,0),\eqno(18)$$
$$(1-(1-\alpha)(1-\beta))P(0,2)=p_2 P(2,0),\eqno(19)$$
$$p_1P(1,0)=\alpha a_1P(0,0)+\alpha a_1\beta P(0,1)+\alpha a_1\beta P(0,2)+\beta P(1,1)+\beta P(1,2),\eqno(20)$$
$$\beta P(1,1)=\alpha a_1(1-\beta)P(0,1),\eqno(21)$$
$$\beta P(1,2)=\alpha a_1(1-\beta)P(0,2),\eqno(22)$$
$$p_2P(2,0)=\alpha a_2P(0,0)+\alpha a_2\beta P(0,1)+\alpha a_2 P(0,2)+
\beta P(2,1)+\beta P(2,2),\eqno(23)$$
$$\beta P(2,1)=\alpha a_2(1-\beta)P(0,1),\eqno(24)$$
$$\beta P(2,2)=\alpha a_2(1-\beta)P(0,2),\eqno(25)$$
$$\sum\limits_{x_1=0}^2\sum\limits_{x_2=0}^2P(x_1,x_2)=1.\eqno(26)$$

Substituting the values expresses  by (3)--(5) for  $P(0,0),$ $P(0,1),$ $P(0,2)$ in (17), we obtain
$$\alpha P^*(0,0)=(1-\alpha)\beta a_1 P^*(0,1)+(1-\alpha)\beta a_2 P^*(0,1).\eqno(27)$$
Taking into account that $a_1+a_2=1,$ we see that (27) is equivalent to (12). Therefore the solution 
(3)--(11) satisfies (17).

Substituting (4), (6) in (18), we obtain
$$(1-(1-\alpha)(1-\beta))a_1P^*(0,1)=
a_1p^*P^*(1,0).\eqno(28)$$
Taking into account (13), we see that (28) is satisfied and therefore the solution (3)--(11) satisfies (18).

Substituting (5), (9) in (19), we obtain
$$(1-(1-\alpha)(1-\beta))a_2P^*(0,1)=
a_2p^*P^*(1,0).\eqno(29)$$
Taking into account (13), we see that (29) is satisfied and therefore the solution (3)--(11) satisfies (19).

Substituting (3)--(8) in (20), we obtain
$$a_1p^*P^*(1,0)=\alpha a_1P^*(0,0)+\alpha a_1^2p^*\beta P^*(0,1)+\alpha a_1a_2p^*\beta P^*(0,1)+\beta a_1^2 p^*(1,1)+\beta a_1a_2 p^*(1,1).\eqno(30)$$
Since the equations $a_1+a_2=1$ and (14) hold, we see that (30) is satisfied and therefore the solution (3)--(11) satisfies (20).

Substitututing (4), (7) in (21), we obtain
$$\beta a_1^2 P^*(1,1)=\alpha a_1^2(1-\beta)P^*(0,1).\eqno(31)$$
From (15) it follows that (31) is satisfied and hence the solution (3)--(11) satisfies (21).

Substituting (5), (8) into (22), we obtain
$$\beta a_1a_2 P^*(1,1)=\alpha a_1 a_2(1-\beta)P^*(0,1).\eqno(32)$$
From (15) it follows that (32) is satisfied and therefore the solution (3)--(11) satisfies (22).

Substituting (3)--(5), (9)--(11) in (23), we obtain
$$a_2p^*P^*(1,0)=\alpha a_2P^*(0,0)+\alpha a_1a_2p^*\beta P^*(0,1)+\alpha a_2^2p^*\beta P^*(0,1)+\beta a_1a_2 P^*(1,1)+\beta a_2^2 P^*(1,1).\eqno(33)$$
From the equations $a_1+a_2=1$ and (13) it follows that (33) is satisfied and therefore the solution (3)--(11) satisfies (23).

Substituting (4), (10) in (24), we obtain
$$\beta a_1a_2 P^*(1,1)=\alpha a_1a_2 P^*(0,1).\eqno(34)$$
From (15) it follows that (34) is satisfied and solution (3)--(11) satisfies (24).

Substituting (5), (11) in (25), we get
$$\beta a_2^2 P^*(1,1)=\alpha a_2^2(1-\beta)P^*(0,1).\eqno(35)$$
From (15) it follows that (35) is satisfied and therefore solution (3)--(11) satisfies (25).

Substituting the solution (3)--(11) in (26) and taking into account that 
$a_1+a_2=1,$ we see that (26) is equivalent to (16).

Thus solution (3)--(11) satisfies all equations of system (17)--(26). Theorem~2 has been proved. 

\vskip 3pt
{\bf Theorem 3.} {\it Suppose the system $S$ contains two cells $(N=2),$ the number of particle types is equal to two $(K=2)$ and the condition $\beta_1=\beta_2=\beta$ is satisfied. Then the sum of the stationary probabilities of the  system $S$ belonging to the set $G(x_1^*,x_2^*)$ is equal to the stationary probability that the system $S^*$ is in the state $(x_1^*,x_2^*),$ $x_1^*=0,1,$ $x_2^*=0,1.$
\vskip 3pt
Proof.} The set $G(0,0)$ contains only the state $(0,0)$ of the system $S.$ From this and (3), the statement of the theorem for the set $G(0,0)$ follows.

The set $G(0,1)$ consists of the states $(0,1)$ and $(0,2)$ of system $S.$
Adding (4) and (5), and taking into account that $a_1+a_2=1,$  we obtain
$$P(0,1)+P(0,2)=a_1P^*(0,1)+a_2P^*(0,1)=P^*(0,1),$$
and hence the stationary probability of system $S$ being in the states of set $G(0,1)$ is equal to the stationary probability of system $S^*$ being in state $S.$

Set $G(1,0)$ consists of states $(1,0)$ and $(2,0)$ of system $S.$
Adding (6), (9), we obtain
$$P(1,0)+P(2,0)=\frac{a_1p^*}{p_1}P^*(0,1)+
\frac{a_2p^*}{p_2}P^*(0,1). \eqno(36)$$
Combining (1), (36), we get the statement of Theorem 3 for the set $G(0,1).$

The set $G(1,1)$ consists of the states $(1,1),$ $(1,2),$ $(2,1)$ and $(2,2)$ of the system $S.$
Adding (7), (8), (10) and (11), we obtain
$$P(1,1)+P(1,2)+P(2,1)+P(2,2)=a_1^2P^*(1,1)+
a_1a_2P^*(1,1)+a_1a_2P^*(1,1)+a_2^2P^*(1,1).\eqno(37)$$
From (37) and the equation $a_1+a_2=1$ it follows that the stationary probability that the system $S$ being in the states of the set $G(1,1)$ is equal to the stationary probability that the system $S^*.$ is in the state $(1,1).$ Theorem~3 has been proved.
\vskip 3pt
Denote by $\rho_i$ $\rho_i^*$ the stationary probability of occupancy (particle density) for cell~$i,$ $i=1,2,\dots,N,$ in systems $S$ and $S^*$ respectively.
\vskip 3pt
{\bf Theorem 4.} {\it The particle densities in the cells 1 and the cell 2 are the same for systems $S$ and
$S^*,$ i.~e. $\rho_i=\rho_i^*,$ $i=1,2.$
\vskip 3pt
Proof.} The value of $\rho_1$ is equal to the probability that the system $S$ is in the state of the set $G(1,0)$ or the set $G(1,1),$ and the value of $\rho_1^*$ is equal to the probability that the system $S^*$ is in the state $(1,0)$ or the state$(1,1).$ From this and Theorem~2, we get Theorem 4.
\vskip 3pt
Denote by $J$ the rate of the particle flow entering the system $S,$ i.~e. the stationary probability that, at the current moment, a particle enters the system. This probability is also equal to the probability that, at the current moment, a particle leaves the system, as well as the probability that, at the current moment, a particle crosses a fixed section of the lattice. We denote by $J^*$ the rate of the particle flow for the system $S^*.$
\vskip 3pt
{\bf Theorem 5.} {\it The particle flow rate the same for the systems $S$ and $S^*,$ i.~e. $J=J^*,$ $i=1,2.$
\vskip 3pt
Proof.} The particle flow rate is equal to the product of the value $\alpha$ and the stationary probability that the cell~1 is vacant, 
and therefore, $J=\alpha(1-\rho),$ $J^*=\alpha(1-\rho^*).$ Taking into account Theorem~4, we get Theorem~5.

\section*{7. Approximate approach}

\hskip 16pt \hskip 16pt Let us propose an approximate approach to evaluate the characteristics of the considered system in the general case. Suppose $N$ is arbitrary. As in the previous sections, we associate the system $S$ under consideration with a system $S^*$ with one type of particles. For the systtem $S^*$ the probability that a particle tries to move at the current moment and the probability that a particle tries to leave the system from the utmost right cell are calculated according to the formulas (1), (2), i.~e. these probabilities are equal to the weighted harmonic average of the corresponding values for different types of particles in the $S$ system.

{\it The approximate values of the densities $\rho_1,\dots,\rho_N$
in the cells of the system $S$ and the rate $J$ in the $S$ system are assumed to be equal, respectively, to the values of the densities $\rho_1^*,\dots,\rho_N^*$ in the cells of the system $S^*$ and rate of the system $S^*$ and the rate $J^*$ in the system $S^*.$
}

Table 1 contains the approximate values computedaccording to the approximate approach
and the exact values (up to rounding) of the evaluated characteristics found from the system
of equations for stationary state probabilities. The number of cells was assumed to be equal to two. Unlike the special case considered in Section 6, the values of $\beta_1$ and $\beta_2$ differ from each other. 
The upper parts of the cells contain exact values, and the lower parts contain approximate values.

\newpage

\vskip 20pt
{\bf Table 1.} Exact and approximate values of density and rate. 
\vskip 10pt
\begin{tabular}{|c|c|c|c|c|c|c|c|c|c|c|}
\hline
id&$\alpha$&$a_1$&$a_2$&$p_1$&$p_2$&$\beta_1$&$\beta_2$&
$\rho_1$&$\rho_2$&$J$\\
\hline
\hline
1&$\frac{2}{5}$&$\frac{3}{7}$&$\frac{4}{7}$&$\frac{3}{5}$&$\frac{4}{5}$&$\frac{3}{10}$&$\frac{2}{5}$&$\frac{0.5149}{0.5142}$&$\frac{0.5544}{0.5552}$&
$\frac{0.1940}{0.1943}$\\
\hline
\hline
2&$\frac{1}{5}$&$\frac{2}{5}$&$\frac{3}{5}$&$\frac{2}{5}$&$\frac{3}{5}$&$\frac{1}{5}$&$\frac{3}{10}$&$\frac{0.4135}{0.4118}$&$\frac{0.4692}{0.4706}$&$\frac{0.1173}{0.1176}$\\
\hline
\hline
3&$\frac{1}{5}$&$\frac{1}{3}$&$\frac{2}{3}$&$\frac{2}{5}$&$\frac{4}{5}$&$\frac{1}{5}$&$\frac{2}{5}$&$\frac{0.3583}{0.3529}$&$\frac{0.4278}{0.4314}$&$\frac{0.1283}{0.1294}$\\
\hline
\hline
4&$\frac{8}{25}$&$\frac{3}{4}$&$\frac{1}{4}$&$\frac{12}{25}$&$\frac{18}{25}$&$\frac{9}{25}$&$\frac{11}{25}$&$\frac{0.4752}{0.4749}$&$\frac{0.4393}{0.4455}$&$\frac{0.1679}{0.1680}$\\
\hline
\hline
5&$\frac{8}{25}$&$\frac{3}{4}$&$\frac{1}{4}$&$\frac{12}{25}$&$\frac{18}{25}$&$\frac{1}{25}$&$\frac{11}{25}$&$\frac{0.5744}{0.5723}$&$\frac{0.5958}{0.6048}$&$\frac{0.1362}{0.1369}$\\
\hline
\end{tabular}
\vskip 10pt
 The proposed approach gives the following results for the case of three cells.
 \vskip 5pt 
Suppose $N=3,$ $\alpha=\frac{1}{5},$ $a_1=\frac{2}{5},$ $a_2=\frac{3}{5},$ $p_1=\frac{2}{5},$ $p_2=\frac{3}{5},$ 
 $\beta_1=\frac{1}{5},$ $\beta_2=\frac{3}{10}.$
The exact values are
$$\rho_1=0.3988,\ \rho_2=0.4374,\ \rho_3=0.4764,\  J=0.1202.$$
The approximate values are
$$\rho_1^*=0.4012,\ \rho_2^*=0.4415,\ \rho_3^*=0.4838,\ J^*=0.1198.$$

\section*{8.Conclusion}

\hskip 18pt

A traffic model on an open one-dimensional lattice with discrete time is considered.
In the system under consideration, there are several types of particles, characterized by the probability of particle attempts to move 
forward at the current moment or to leave the system. For a particular case, theorems on stationary probabilities of system states, 
particle density in a cell depending on its location and particle flow intensity are proved. For the general case of the system, an 
approximate approach is proposed. According to this approach, the density values in cells and the values of rate are assumed to be equal to 
the corresponding values for a system with one type of particles. In this system, the probabilities of attempts to move or leave the system
are assumed to be equal to the average harmonic values of the corresponding probabilities of such attempts for particles of different types in the original system.

\bibliographystyle{unsrt}

\begin{thebibliography}{1}

\bibitem{Nagel}
Nagel K., Schreckenberg M. A cellular automation models
for freeway traffic. J.~Phys. I. France 2, 1992, pp. 2221--2229. DOI:~10.1051/jp1:1992277

\bibitem{Wolfram}
Wolfram S. Statistical mechanics of cellular automata. Reviews of Modern Physics, 55, 3, 601, 1983. DOI:~10.1103/RevModPhys.55.601

\bibitem{Spitzer}
Spitzer F. Interaction of Markov procsses //In Random Walks, Browing Motion, and Interacting Particle Systems: A Festschrift in Honor of Frank Spitzer. – Boston, MA : Birkhäuser Boston, 1991, pp. 66-110.

\bibitem{Blank}
Blank M.L. Metric properties of discrete time exclusion type processes in continuum. Journal of Statistical Physics, 2010, vol.~140, no.~1, pp.~170--197. 

\bibitem{Schadschneider}
Schadschneider A., Schreckenberg M. Cellular automation models and traffic flow //Journal of Physics A: Mathematical and General, 1993, vol.~26, no.~15, pp.~L679.

\bibitem{Kanai}
Kanai~M., Nishinari~K., Tokihiro~T. (2006). Exact solution and asymptotic behaviour of the asymmetric simple exclusion process on a ring. Journal of Physics A: Mathematical and General, 39(29), 9071.
DOI: 10.1088/0305-4470/39/29/004

\bibitem{Buslaev}
Buslaev A.P., Tatashev A.G. Particles flow on the regular polygon //Journal of
Concrete and Applicable Mathematics, 2011, vol.~9, no.~4, pp.~290--303.

\bibitem{Kanai}
Kanai M. Exact solution of the zero-range process: fundamental diagram of the corresponding exclusion process //Journal of Physics A: Mathematical and Theoretical, 2007,  vol.~40,  no.~26,  pp.~7127. DOI: 10.1088/1751-8113/40/26/001

\bibitem{Tatashev}
 Kozlov V.V., Buslaev A.P., Tatashev A.G., Yashina M.V.
Dynamical systems on 
honeycombs, Traffic and granular flow '13, Springer, 2015, pp. 441–452 

\bibitem{Belitsky}
Belitsky V., Ferrari P. A. Invariant measures and convergence properties for cellular automaton 184 and related processes. Journal of Statistical Physics. 2005, 118,  pp.~589--623.

\bibitem{Bugaev}
Bugaev A.S., Tatashev A.G., Yashina M.V., Lavrov O.S., Nosov E.A. (2019). Reconstruction of traffic flow dynamics based on deterministicstochastic model and data obtained from intelligent transport systems. T-Comm, 13(10), 35-44.

\bibitem{Daduna}
Daduna H. Queueing networks with discrete time scale: explicit expressions for the steady behavior of discrete time stochastic networks. Vol.~2046. Springer, 2003. 

\bibitem{Gordon}
Gordon W.J., Newell G.F. Closed queuing systems with exponential servers //Operations research, 1967, vol.~15,  no~2, pp.~254--265. DOI: 10.1287/opre.15.2.254

\bibitem{Kelly}
Kelly F.P. Networks of queues //Advances in Applied Probability, 1976, vol.~8, no.~2, pp. 416--432. DOI: 10.2307/1425912

\bibitem{Yashina}
Yashina M., Tatashev A. Traffic model based on synchronous and asynchronous exclusion processes //Mathematical Methods in the Applied Sciences,  2020,  vol.~43,  no.~14,  pp.~8136--8146. DOI: 10.1002/mma.6237

\bibitem{Yashina}
A.Tatashev, M. Yashina. Self-organization of two-contours dynamical system with common node and cross movement. WSEAS Trans. Math. 18(45) (2019) 373.

\bibitem{Evans}
Evans M.R., Rajewsky N., Speer E.R. Exact solution of a cellular automaton for traffic //Journal of Statistical Physics,1999, vol. 95, pp.~45--96.

\bibitem{Derrida}
Derrida B. An exactly soluble non-equilibrium system: The asymmetric simple exclusion process. Physics Reports, 1998, vol.~301, no.~1--3, pp.~65--83. \\
DOI: 10.1016/S0370-1573(98)00006-4 .

\bibitem{Xiao}
Xiao S., Liu M., Cai J.J. (2009). Asymmetric coupling in two-lane simple exclusion processes: Effect of unequal injection rates. Physics Letters A, 374(1), 8-12.

\bibitem{Kemeny}
Kemeny J.G., Snell J.L. Finite Markov Chains. Springer-Verlag, New York, Berlin, Heidelberg Tokyo, 1975. 



\end{thebibliography}

\end{document}